\documentclass[11pt,reqno,oneside]{amsart}
 \usepackage{latexsym}
 \usepackage{amssymb}
 \usepackage{amsfonts}
\usepackage{graphics}
\usepackage{mathcomp}
 \author[Mathonet]{P. Mathonet}
\thanks{University of Li\`ege, Institute of mathematics, Grande Traverse, 12 - B37, B-4000
 Li\`ege, Belgium email : P.Mathonet@ulg.ac.be, Fabian.Radoux@ulg.ac.be}
\author[Radoux]{F. Radoux}
\date{\today} 
\title[On natural and invariant quantizations]{Existence of natural and conformally invariant quantizations of arbitrary symbols}

\newtheorem{lem}{Lemma}
\newtheorem{thm}[lem]{Theorem}
\newtheorem{prop}[lem]{Proposition}

\theoremstyle{remark}
\newtheorem{rem}{Remark}
\theoremstyle{definition}
\newtheorem{defi}{Definition}

\newcommand{\R}{\mathbb{R}}
\newcommand{\Z}{\mathbb{Z}}

\newcommand{\N}{\mathbb{N}}

\newcommand{\V}{\mathcal{V}}
\newcommand{\D}{\mathcal{D}}
\renewcommand{\S}{\mathcal{S}}

\newcommand{\g}{\mathfrak{g}}

\newcommand{\h}{\mathfrak{h}}

\newcommand{\euler}{\mathcal{E}}

\newcommand{\tree}[1][\gamma]{\mathcal{T}_{{#1}}}

\newcommand{\gm}{{\g_{-1}}}
\newcommand{\Cf}{C^\flat}

\newcommand{\DM}{{\mathcal D}(M)}
\newcommand{\DkM}{{\mathcal D}^k(M)}
\newcommand{\SkM}{{\mathcal S}^k(M)}
\newcommand{\SM}{{\mathcal S}(M)}
\begin{document}
\begin{abstract}
A quantization over a manifold can be seen as a way to construct a differential operator with prescribed principal symbol. The map from the space of principal symbols to the space of differential operators is moreover required to be a linear bijection. 

It is known that there is in general no natural quantization procedure. 
However, considering manifolds endowed with additional structures, such as projective or pseudo-conformal structures, one can seek for quantizations that depend on this additional structure and that are natural if the dependence with respect to the structure is taken into account. The existence of such a quantization was conjectured by P. Lecomte in \cite{Leconj} in the context of projective and conformal geometry.

The question of existence of such a quantization was addressed in a series of papers in the context of projective geometry, using the framework of Thomas-Whitehead connections (see for instance \cite{Bor,Hansoul,Fox,sarah}). In \cite{MR,MR1}, we recovered the existence of a quantization that depends on a projective structure and that is natural (provided some critical situations are avoided), using the theory of Cartan projective connections. 

In the present work, we show that our method can be adapted to pseudo-conformal geometry to yield the so-called natural and conformally invariant quantization for arbitrary symbols, still outside some critical situations. 

Moreover, we give new and more general proofs of some results of \cite{MR1} and eventually, we notice that the method is general enough to analyze the problem of natural and invariant quantizations in the context of manifolds endowed with irreducible parabolic geometries studied in \cite{CapGov}. 
\end{abstract}
%%%
\maketitle
%%%
\section{Introduction}
%%%
%%%
Consider a manifold $M$ and a space $\DM$ of differential operators acting between spaces of tensor fields over $M$.
The space $\DM$ is filtered by the order of differential operators :
\[\DM=\cup_{k=0}^\infty\DkM,\]
where $\DkM$ is the space of operators of order at most $k$. 

The principal symbol of a differential operator of order $k$ is a well defined tensor field : there is a natural map --the principal symbol operator--
\[\sigma_k:\DkM\to \SkM,\]
and the space of symbols associated to $\DM$ is then the graded space
\[\SM=\oplus_{k=0}^\infty\SkM.\]
The purpose of the so-called ``natural quantization problem'' is to construct a differential operator with prescribed symbol in a natural way. In this framework, a natural quantization is a natural bijection $Q$ from $\SM$ to $\DM$ such that $\sigma_k\circ Q=\mathrm{Id}$ on $\SkM$ for every $k$.

It is known that in general, such a natural quantization does not exist. There are two possible ways to weaken the naturality condition. The first one leads to the concept of $G$-equivariant quantizations on manifolds endowed with the action of a Lie group $G$. The second one is to consider the more general notion of natural and invariant quantizations.

The concept of $G$-equivariant quantization was defined by P. Lecomte and V. Ovsienko in \cite{LO} in the following way : if a Lie group $G$ acts (locally) on a manifold $M$, the action can be lifted to tensor fields and to differential operators and symbols. A $G$-equivariant quantization is then a quantization that exchanges the actions of $G$ on symbols and differential operators.
In \cite{LO}, the authors considered a space of differential operators ${\mathcal D}_\lambda(\R^m)$ acting on $\lambda$-densities and the projective group $G=PGL(m+1,\R)$ acting on $\R^m$ by linear fractional transformations. They showed that there exists a unique
\emph{projectively equivariant quantization}.

In \cite{DO}, the authors studied the spaces $\mathcal{D}_{\lambda \mu}(\R^m)$ of
differential operators transforming $\lambda$-densities into $\mu$-densities.
They showed the existence and uniqueness of a projectively equivariant
quantization, provided the shift value $\delta=\mu-\lambda$ does not belong to
a set of critical values, and in \cite{BHMP}, a first example of projectively equivariant quantizations for differential
operators acting on tensor fields was considered.

In \cite{DLO}, the authors considered the group $SO(p+1,q+1)$ acting on the
space $\R^{p+q}$ or on a manifold endowed with a flat pseudo-conformal structure of signature $(p,q)$.
They also showed the existence and uniqueness of a \emph{conformally
equivariant quantization} provided the shift value is not critical.

The algebras of fundamental vector fields associated to the actions of $PGL(m+1,\R)$ and of $SO(p+1,q+1)$ turn out to be maximal in the set of proper subalgebras of polynomial vector fields. In \cite{BMMax}, the authors classified the subalgebras having this property. They correspond to simple algebras $\g$ carrying a $|1|$-grading $\g=\gm\oplus\g_0\oplus\g_1$. Finally, in \cite{BM}, the methods of \cite{DLO} were adapted to obtain existence of the quantization for most of these algebras $\g$ (and the corresponding groups $G$) provided the shift value is not critical.

The concept of natural and invariant quantization appeared in the conformal case, in \cite{DO1} and \cite{Loubon}, where the authors showed that the conformally equivariant quantization procedure for symbols of degree two and three can be expressed using the Levi Civita connection associated to a pseudo-Riemannian metric in such a way that it only depends on the \emph{conformal class} of the metric.

In the projective case, in \cite{Bou1,Bou2}, S. Bouarroudj showed 
that the formula for the projectively equivariant quantization for differential operators of order two
 and three could be expressed using a torsion-free linear connection, in such
 a way that it only depends on the \emph{projective
 class} of the connection.
 
In \cite{Leconj}, P. Lecomte conjectured the existence of a quantization 
procedure  $Q:\S(M)\to{\mathcal{D}}_{\lambda,\mu}(M)$
depending on a torsion-free linear connection (resp. pseudo-Riemannian metric), that would be
 natural in all arguments and that would remain invariant under a projective (resp. conformal)
change of connection (resp. metric).

The existence of such a \emph{Natural and projectively invariant
   quantization procedure} was first proved by M. Bordemann 
in \cite{Bor}, using the notion of Thomas-Whitehead connection  associated to
   a projective class of connections.
This result was generalized in a series of recent papers. First, 
 S. Hansoul adapted the construction in \cite{Hansoul} in order to extend the results of \cite{BHMP}.
Then, in \cite{MR}, we gave an alternative proof of the results of M. Bordemann, using the theory of projective Cartan connections. 

In \cite{Fox,sarah}, these results were generalized to deal with multilinear operators or linear operators acting on arbitrary tensors.

At the same time, in \cite{MR1}, we gave a proof of the existence of natural and projectively invariant quantizations using Cartan projective connections. One of the advantages of our method was to obtain a formula showing a close relationship with the formulae obtained in the context of $PGL(m+1,\R)$-equivariant quantizations over $\R^m$.

In the conformal situation, we proved in \cite{MR2} the existence --outside critical situations-- of the conformally invariant quantization for symbols at most four and for differential operators acting between densities.  

In this paper, we give a general result of existence of a conformally invariant natural quantization for differential operators acting on arbitrary tensors and for any order of differentiation, provided the situation is not critical. We actually adapt the tools of \cite{MR1} to the conformal case. We also give new and more general proofs of some of the results of \cite{MR1}, and it turns out that our arguments and tools are general enough to obtain a natural and invariant quantization associated to most of geometric structures corresponding to simple $|1|$-graded algebras.

%The paper is organized as follows. We define the general concept of conformally
%invariant and natural quantizations over arbitrary manifolds in section \ref{tens}. 
%In section \ref{curvtools}, we adapt the tools of \cite{MR1} to the conformal situation, and
%we eventually use these tools in order to build the quantization in section \ref{curvconst}. 
%%%%
%%%%
\section{Problem setting}\label{tens}
In this section, we will describe the definitions of the spaces of
differential operators acting on tensor fields and of their corresponding
spaces of symbols. Then we will
set the problem of existence of natural and conformally invariant
quantizations. Throughout this work, we let $M$ be a smooth manifold of dimension $m\geq 3$.
\subsection{Tensor fields}
The arguments of the differential operators that we will 
consider are classical tensor fields. Recall that one may see them as sections of
vector bundles associated to the linear frame bundle $P^1M$.
We consider irreducible representations of the group $GL(m,\R)$ defined as
follows~:
let $(V,\alpha_D)$ be the representation of $GL(m,\R)$ corresponding to a Young
diagram $Y_D$ of depth $n<m$. Fix $\lambda\in\R$ and $z\in \Z$ and set 
\begin{equation}\label{rep}\alpha(A)u=\vert det(A)\vert^{\lambda} 
(det(A))^z\alpha_D(A)u,\end{equation}
for all $A\in GL(m,\R)$, and $ u\in V$.\\
If $(V,\alpha)$ is such a representation, we denote by $V(M)$ the vector bundle
\[P^1M\times_{\alpha}V.\]
We denote by $\V(M)$ the space of smooth sections of $V(M)$. 
 This space can be identified with the space $C^{\infty}(P^1M,V)_{GL(m,\R)}$ of $GL(m,\R)$-invariant smooth functions, i.e.
functions $f$ such that
\[f(u A) = \alpha(A^{-1}) f(u)\quad \forall u \in P^1M,\;\forall A\in
GL(m,\R).\]
Finally, since $V(M)\to M$ is  associated to $P^1M$, the group
$\mathrm{Diff}(M)$ of diffeomorphisms of $M$ acts in a well-known manner on the space $\V(M)$. 
%%%
\subsection{Differential operators and symbols}\label{Diffops}
%%%%
If $(V_1,\rho_1)$ and $(V_2,\rho_2)$ are representations of $GL(m,\R)$, we
 denote by $\mathcal{D}(\V_1(M),\V_2(M))$ (or simply by $\D(M)$ if there is no
 risk of confusion) the space of 
linear differential operators from $\V_1(M)$ to
$\V_2(M)$. 
The actions of $\mathrm{Vect}(M)$ and
$\mathrm{Diff}(M)$ are induced by their actions on $\V_1(M)$ and
 $\V_2(M)$. For instance, one has
\[(\phi\cdot D)(f) = \phi\cdot(D(\phi^{-1}\cdot f)),\quad\forall f\in \V_1(M)
, D\in \mathcal{D}(M),\mbox{and}\, \phi
\in \mathrm{Diff}(M).\]
The space $\mathcal{D}(M)$ is filtered by the order of
differential operators. We denote by $\mathcal{D}^k(M)$ the
space of differential operators of order at most $k$. 
The space of \emph{symbols}, which we will denote by $\S_{V_1,V_2}(M)$ or
simply by $\S(M)$, is then the graded space associated to $\mathcal{D}(M)$.

We denote by $S^l_{V_1,V_2}$ the vector space
 $S^l\R^m\otimes gl(V_1,V_2)$. We denote by $\rho$ the natural representation of $GL(m,\R)$
  on this space (the representation of $GL(m,\R)$ on
 symmetric tensors is the natural one). We then denote by 
$S^l_{V_1,V_2}(M)\to M$
 the vector bundle
\[P^1M\times_{\rho}S^l_{V_1,V_2}\to M,\]
and by $\mathcal{S}^l_{V_1,V_2}(M)$ the space of smooth sections of
$S^l_{V_1,V_2}(M)\to M$, that is, the space 
$C^{\infty}(P^1M,S^l_{V_1,V_2})_{ GL(m,\R)}$.

 Then the \emph{principal symbol operator} $\sigma_l :
\mathcal{D}^l(M)\to \mathcal{S}^l_{V_1,V_2}(M)$ commutes with
the action of diffeomorphisms and is a bijection from the quotient space
$\mathcal{D}^l(M)/\mathcal{D}^{l-1}(M)$ to
$\mathcal{S}^l_{V_1,V_2}(M)$.
Hence the space of symbols is nothing but
\[\mathcal{S}(M)
=\bigoplus_{l=0}^{\infty}\mathcal{S}^l_{V_1,V_2}(M),\]
endowed with the classical action of $\mathrm{Diff}(M)$.
\subsection{Natural and invariant quantizations}
A \emph{quantization on $M$} is a linear bijection $Q_{M}$
  from the space of symbols $\mathcal{S}(M)$ to the space of differential operators $\mathcal{D}(M)$ such
  that
\[\sigma_k(Q_{M}(S))=S,\;\forall S\in\mathcal{S}^{k}_{V_1,V_2}(M),\;\forall
  k\in\mathbb{N}.\]
In the conformal sense, a \emph{natural quantization} is a collection of quantizations $Q_M$ depending on a pseudo-Riemannian metric  such that 
\begin{itemize}
\item For all pseudo-Riemannian metric $g$ on $M$, $Q_{M}(g)$ is a quantization,
\item If $\phi$ is a local diffeomorphism from $M$ to $N$, then one has
\[Q_{M}(\phi^{*}g)(\phi^{*}S)=\phi^{*}(Q_{N}(g)(S)),\]
 for all pseudo-Riemannian metrics $g$ on $N$, and all
$S\in\mathcal{S}(N).$
\end{itemize}
Recall now that two pseudo-Riemannian metrics $g$ and $g'$ on a manifold $M$ are conformally equivalent if and only if there exists a positive function $f$ such that $g'=fg$. 

A quantization $Q_M$ is then \emph{conformally invariant} if one has $Q_M(g)=Q_M(g')$ whenever $g$ and $g'$ are conformally equivalent.
\section{Conformal group and conformal algebra}
Let us now recall the definition of the algebraic objects that we will use throughout this work. The grading of these objects are of special importance.
\subsection{The conformal group}\label{projconf}
%%%%%
Given $p$ and $q$ such that $p+q=m$, we consider the bilinear symmetric form of signature $(p+1,q+1)$ on $\R^{m+2}$ defined by
\[B:\R^{m+2}\times\R^{m+2}\to\R:(x,y)\mapsto {^ty}Sx,\]
where $S$ is the matrix of order $m+2$ given by
\[S=
      \left(
        \begin{array}{ccc}
           0 & 0 & -1\\
           0 & J & 0\\
          -1 & 0 & 0
       \end{array}
     \right),
    \]
and 
\[J=
      \left(
        \begin{array}{cc}
           I_{p} & 0\\
           0     & -I_{q}
          \end{array}
     \right)
    \]
    represents a nondegenerate symmetric bilinear form $g_0$ on $\R^m$, namely
    \[g_0:\R^{m}\times\R^{m}\to\R:(x,y)\mapsto {^ty}Jx.\]    
As we continue, we will use the classical isomorphism between $\R^m$ and $\R^{m*}$ defined by the symmetric bilinear form represented by $J$ :
\[\sharp : \R^m\to \R^{m*}: x\mapsto x^\sharp : x^\sharp(y)= g_0(x,y)\]
and $\flat=\sharp^{-1}$, and we also denote by $|x|^2$ the number $g_0(x,x)$.

The M\"{o}bius space is the projection of the light cone associated to the metric $B$ on the projective space $\mathbb{R}P^{m+1}$.

We consider the group $G$ made of linear transformations that leave $B$ invariant, modulo its center, that is,
\[G=\{X\in GL(m+2,\mathbb{R}) : {^tX}SX=S\}/\{\pm I_{m+2}\}.\]
It acts transitively on the M\"{o}bius space $S^{m}$.

The group $H$ is the
isotropy subgroup of $G$ at the point $[e_{m+2}]$ of the M\"{o}bius space :
\[H=\{ \left( 
           \begin{array}{ccc}
              a^{-1} & 0 & 0\\
              a^{-1}A\xi^\flat & A & 0\\
              \frac{1}{2a}|\xi|^2 & \xi & a
           \end{array}
        \right): A\in
        O(p,q),a\in\mathbb{R}_{0},\xi\in\mathbb{R}^{m*}\}/\{\pm I_{m+2}\}.\]
As in the projective situation, $H$ is a semi-direct product $G_{0}\rtimes
        G_{1}$. Here $G_{0}$ is isomorphic to $CO(p,q)$ and $G_{1}$ is
        isomorphic to $\mathbb{R}^{m*}$. There is also a projection 
\[\pi:H\mapsto CO(p,q):\left[
        \left( 
           \begin{array}{ccc}
              a^{-1} & 0 & 0\\
             a^{-1}A\xi^\flat & A & 0\\
               \frac{1}{2a}|\xi|^2 & \xi & a
           \end{array}
        \right)\right]\mapsto \frac {A}{a}.\]
        \subsection{The conformal algebra}\label{algebra}
The Lie algebra of $G$ is $\mathfrak{g}=so(p+1,q+1)$. It decomposes as a direct sum of subalgebras :
\begin{equation}\label{grading}\g=\mathfrak{g}_{-1}\oplus\mathfrak{g}_{0}\oplus\mathfrak{g}_{1}\end{equation}
where $\g_{-1}\cong\R^{m}$, $\g_0\cong co(p,q)$, and $\g_1\cong\mathbb{R}^{m*}$.
The isomorphism is given explicitly by
\[
         \left( 
           \begin{array}{ccc}
              -a & v^\sharp & 0\\
              \xi^{\flat} & A & v\\
              0 & \xi & a
           \end{array}
        \right)\mapsto(v,A-aI_{m},\xi).\]
This correspondence induces a structure of Lie algebra on $\mathbb{R}^{m}\oplus
co(p,q)\oplus\mathbb{R}^{m*}$. It is easy to see that the adjoint actions $G_0$ and of $co(p,q)$ on $\gm=\R^m$ and on $\g_1={\R^m}^*$ coincides with the natural actions of $CO(p,q)$ and of $co(p,q)$. 
%Moreover, one has
%\begin{equation}\label{bra} [h,x]=-x\otimes h + h^\flat\otimes x^\sharp -\langle h,x\rangle Id.\end{equation}
%
The Lie algebras corresponding to $G_{0}$,
$G_{1}$ and $H$ are respectively $\mathfrak{g}_{0}$, $\mathfrak{g}_{1}$, and
$\mathfrak{g}_{0}\oplus\mathfrak{g}_{1}$.

Actually, the simple Lie algebras carrying a grading (\ref{grading}) are known as simple $|1|$-graded algebras or Irreducible Filtered Lie algebras of Finite Type (IFFT algebras for short). Hence the algebra $so(p+1,q+1)$ is a particular case of such an algebra. The classification of these algebras was obtain in \cite{Koba1}. Recall also that for every simple $|1|$-graded algebra, the subalgebra $\g_0$ is reductive and decomposes as 
\begin{equation}\label{euler}
\g_0=\h_0\oplus\R \euler\end{equation}
where $\h_0$ is semi-simple and where the \emph{grading} or \emph{Euler} element $\euler$ is defined by $ad(\euler)\vert_{\g_{k}}=k Id$ ($k\in\{-1,0,1\}$) and is therefore in the centre of $\g_{0}$. It is also noteworthy that $\g_{-1}$ and $\g_1$ are always dual of each other because of the non degeneracy of the Killing form $K$ of $\g$. 

Let us now close this section by to technical results about Killing-dual bases in such algebras. The first one is taken directly from \cite{BM}.
\begin{prop}\label{basis}
For every IFFT-algebra $\g$, there exists a basis $(e_i, \euler,A_j,\varepsilon^i)$ of $\g$ such that $(e_i:1\leq i\leq d)$ is a basis of $\g_{-1}$, $(A_j : 1\leq j\leq \dim \h_0)$ is a basis of $\h_0$ and $(\varepsilon^i:1\leq i\leq d)$ is a basis of $\g_1$ and such that the Killing dual basis writes $(\varepsilon^i,\frac{1}{2d}\euler,A_j^*,e_i)$, where $A_j^*$ is in $\h_0$ for all $j\leq \dim \h_0$.
Moreover we have 
\[\sum_{i=1}^d[\varepsilon^i,e_i]=\frac{1}{2}\euler.\]
\end{prop}
As we continue, we will need other relations concerning these bases.
\begin{prop}\label{rel}
Denote respectively by $a_j$ and $a_j^*$ the matrix representations of $ad(A_j)|_\gm$ and $ad(A_j^*)|_\gm$ in the basis $(e_i:i\leq d)$. Then we have 
\[\begin{array}{lll}[\varepsilon^r,e_i]&=&\frac{1}{2d}\delta_i^r\euler-\sum_{j=1}^{\dim \h_0}a_{ji}^rA_j^*\\
&=&\frac{1}{2d}\delta_i^r\euler-\sum_{j=1}^{\dim \h_0}a_{ji}^{*r}A_j
\end{array}\]
for every $i,r\leq d$.
Moreover we have 
\[[A_j,\varepsilon^r]=-\sum_{k=1}^da_{jk}^r\varepsilon^k\quad\mbox{and}\quad [A_j^*,\varepsilon^r]=-\sum_{k=1}^da_{jk}^{*r}\varepsilon^k.\]
\end{prop}
\begin{proof}
Let us prove the first relation. The other ones are obtained in a similar way. We know that $[\varepsilon^r,e_i]$ belongs to $\g_0$ since the algebra $\g$ is $|1|$-graded. Hence we have a decomposition
\[[\varepsilon^r,e_i]=b^r_i\euler+\sum_{j=1}^{\dim \h_0}b_{ji}^rA_j^*.\]
We compute the coefficients using the Killing dual basis given in Proposition \ref{basis} : we have 
\[b^r_i=K([\varepsilon^r,e_i],\frac{1}{2d}\euler)=-K(\varepsilon^r,\frac{1}{2d}[\euler,e_i])=\frac{1}{2d}K(\varepsilon^r,e_i),\]
by the invariance of the Killing form $K$ and the definition of $\euler$.
In the same way we obtain 
\[b_{ji}^r=K([\varepsilon^r,e_i],A_j)=-K(\varepsilon^r,[A_j,e_i])=-K(\varepsilon^r,\sum_sa_{ji}^se_s)=-a_{ji}^r.\]
\end{proof}
%%%%%
\section{Cartan fiber bundles and connections}\label{curvtools}
\subsection{Cartan fiber bundles}
It is well-known that there is a bijective and natural correspondence between
the conformal structures on $M$ and the reductions of $P^{1}M$ to the
structure group $G_{0}\cong CO(p,q)$. The representations $(V,\alpha)$ of $GL(m,\R)$ defined so far can be restricted to the group $CO(p,q)$. Therefore, once a conformal structure is given, i.e. a reduction $P_0$ of $P^1M$ to $G_0$, we can identify tensors fields of type $V$ as $G_0$ invariant functions on $P_0$.

In \cite{Kobabook}, one shows that it is possible to
associate at each $G_{0}$-structure $P_{0}$ a principal $H$-bundle $P$ on $M$,
this association being natural and obviously conformally invariant. Since $H$ can be considered as a subgroup of $G^2_m$, this $H$-bundle can
be considered as a reduction of $P^{2}M$. The relationship between conformal structures and reductions of $P^2M$ to $H$ is given by the following proposition.
\begin{prop}
There is a natural one-to-one correspondence between the conformal equivalence
classes of pseudo-Riemannian metrics on $M$ and the reductions of $P^2M$ to $H$.\end{prop}
Throughout this work, we will freely identify conformal structures and reductions of $P^2M$ to $H$. 
\subsection{Cartan connections}
Let $L$ be a Lie group and $L_{0}$ a closed subgroup. Denote by $\mathfrak{l}$
and $\mathfrak{l}_{0}$ the corresponding Lie algebras. Let $N\to M$ be a
principal $L_{0}$-bundle over $M$, such that $\dim M=\dim L/L_{0}$. A Cartan
connection on $N$ is an $\mathfrak{l}$-valued one-form $\omega$ on $N$ such
that
\begin{enumerate}
\item
If $R_{a}$ denotes the right action of $a\in L_0$ on $N$, then
$R_{a}^{*}\omega=Ad(a^{-1})\omega$,
\item
If $k^{*}$ is the vertical vector field associated to $k\in\mathfrak{l}_0$,
then $\omega(k^{*})=k,$
\item
$\forall u\in N$, $\omega_{u}:T_{u}N\mapsto\mathfrak{l}$ is a linear
  bijection.
\end{enumerate}
When considering in this definition a principal $H$-bundle $P$, and taking as group $L$ the group $G$ and for $L_0$ the group $H$, we obtain the definition of Cartan conformal connections.

If $\omega$ is a Cartan connection defined on an $H$-principal bundle $P$, then its
curvature $\Omega$ is defined by
 \begin{equation}\label{curv} 
\Omega = d\omega+\frac{1}{2}[\omega,\omega].
\end{equation}
The notion of \emph{Normal} Cartan connection is defined by natural conditions
imposed on the components of the curvature.

Now, the following result (\cite[p. 135]{Kobabook}) gives the relationship between conformal
structures and Cartan connections :
\begin{prop}
 A unique normal
 Cartan conformal connection is
 associated to every conformal structure $P$. This association is natural.
\end{prop}
The connection associated to a conformal structure $P$ is called the normal conformal connection of the conformal structure.

%\bf{Remark} : \rm{Actually}, all the constructions we are going to do in the sequel can be done for all the AHS-structures defined in \cite{}. In fact, we are going to find natural applications that associate a quantization to a reduction of $P^{2}M$ to $H$, where $H$ is a group corresponding to a certain AHS-structure. In the situations different from the projective or the conformal one, we have no geometric interpretation for these problems.
\section{Lift of equivariant functions}\label{Lift}
%%%%\subsection{Lift of equivariant functions}\label{Lift}
In the previous section, we recalled how to associate an $H$-principal bundle $P$ to a conformal structure $P_0$. We now recall how the tensor fields, and in particular symbols, can be regarded as equivariant functions on $P$.

If $(V,\alpha)$ is a representation of $G_0$, then we may extend it to a representation $(V,\alpha')$ of $H$ by
\[\alpha'=\alpha\circ\pi.\]
Now, using the representation $\alpha'$, we can recall the relationship between
equivariant functions on $P_{0}$ and equivariant functions on $P$ (see \cite{Capcart}): if we
denote by $p$ the projection $P\to P_{0}$ , we have
\begin{prop}
If $(V,\alpha)$ is a representation of $G_0$, then the map
$$p^{*}:C^{\infty}(P_{0},V)\mapsto C^{\infty}(P,V):f\mapsto f\circ p$$
defines a bijection from $C^{\infty}(P_{0},V)_{G_0}$ to
$C^{\infty}(P,V)_{H}$.
\end{prop}
Now, since $\g_{-1}\cong\R^m$ and $\g_1\cong\R^{m*}$ are natural representations of $G_0\cong CO(p,q)$,
they become representations of $H$ and we can state an important property of the invariant differentiation :
\begin{prop}\label{gonabla}If $f$ belongs to $C^{\infty}(P,V)_{G_0}$ then
  $\nabla^{\omega}f\in C^{\infty}(P,\R^{m*}\otimes V)_{G_0}$.
\end{prop}
\begin{proof}The result is a direct consequence of the Ad-invariance of the
  Cartan connection $\omega$. 
\end{proof}
The main point that we will discuss in the next sections is that this result
is not true in general for $H$-equivariant functions : for an $H$-equivariant
function $f$, the function
$\nabla^{\omega}f$ is in general not $G_1$-equivariant.

As we continue, we will use the representation $\rho'_*$ of the Lie algebra of
$H$ on $V$. If we recall that this algebra is isomorphic to
 $\g_{0}\oplus\g_{1}$
then we have 
\begin{equation}\label{rho}\rho'_* (A, \xi) = \rho_*(A),\quad\forall A\in\g_{0}, \xi\in \g_{1}.\end{equation}
In our computations, we will make use of the infinitesimal version of the
equivariance relation : If $f\in C^{\infty}(P,V)_H$ then one has
\begin{equation}\label{Invalg}
L_{h^*}f(u) + \rho'_*(h)f(u)=0,\quad\forall h\in\g_{0}\oplus\g_{1},
\forall u\in P.
\end{equation}
\section{The application $Q_{\omega}$}
The construction of the application $Q_{\omega}$ is based
on the concept of invariant differentiation developed in
\cite{Capinv,Capcart}.
Let us recall the definition :

\begin{defi}
If $f\in C^{\infty}(P,V)$ then $(\nabla^{\omega})^k f
\in C^{\infty}(P,\otimes^k\R^{m*}\otimes V)$ is defined by
\[(\nabla^{\omega})^k f(u)(X_1,\ldots,X_k) = L_{\omega^{-1}(X_{1})}\circ\ldots\circ L_{\omega^{-1}(X_{k})}f(u)\]
for $X_1,\ldots,X_k\in\R^m$.
\end{defi}

\begin{defi}The map $Q_{\omega}$ is defined by its
  restrictions to $ C^{\infty}(P,\otimes^k\g_{-1}\otimes gl(V_1,V_2))$, $(k\in \N)$ : we set
\begin{equation}\label{Qom}
Q_{\omega}(T)(f)=\langle T,(\nabla^{\omega})^k f\rangle,
\end{equation}
for all $T\in C^{\infty}(P,\otimes^k\g_{-1}\otimes gl(V_1,V_2))$ and $f\in C^{\infty}(P,V_1)$.
\end{defi}
%%%%%
Explicitly, when the symbol $T$ writes $t A\otimes h_1\otimes\cdots\otimes h_k$ 
for $t\in C^{\infty}(P)$, $A\in V_1^*\otimes V_2$ and $h_1,\cdots,
h_k\in\R^m\cong\g_{-1}$ then one has
\[Q_{\omega}(T)f=t A\circ
L_{\omega^{-1}(h_{1})}\circ\cdots\circ L_{\omega^{-1}(h_{k})}f,\]
where $t$ is considered as a multiplication operator\vspace{0,2cm}.

\begin{rem}\label{remf} If $T\in C^{\infty} 
(P,\otimes^k\g_{-1}\otimes gl(V_1,V_2))$ is
$H-$equivariant, the differential operator $Q_{\omega}(T)$ does not transform
$H-$equivariant functions into $H-$equivariant functions. Indeed, when $f$ is
$H-$equivariant, the function $(\nabla^{\omega})^k f$ is only
$G_0-$equivariant. Hence the function $Q_{\omega}(T)f$ does not correspond to
a section of $\mathcal{V}_2(M)$. As we continue, we will show that one can modify the
symbol $T$ by lower degree correcting terms in order to solve this problem.
\end{rem}
\section{Measuring the default of equivariance}
Throughout this section, $T$ will denote an element of $C^{\infty}
(P,S^k_{V_1,V_2})_{G_0}$ and $f\in C^{\infty}(P,V_1)_{G_0}$ 
(remark that this ensures that
$Q_{\omega}(T)(f)$ is in $C^{\infty}(P,V_2)_{G_0}$). 
Now, in order to analyze the invariance of functions, we have this first easy
result, which follows from the fact that $\g_1$ is a vector space.
\begin{prop}[\cite{MR1}]If $(V,\alpha)$ is a representation of $G_0$ and becomes a
representation of $H$ as stated in section \ref{Lift}, then a function $v\in
 C^{\infty}(P,V)$ is $H-$equivariant iff 
\[\left\{\begin{array}{l}
\mbox{$v$ is $G_0-$equivariant}\\
\mbox{One has $L_{h^*}v=0$ for every $h$ in $\g_1$}
\end{array}\right.\]
\end{prop}
\subsection{The map $\gamma$}
As we continue, we are interested in measuring the failure of equivariance of the map $Q_{\omega}$. To this aim, we compute the Lie derivative of the differential operator $Q_{\omega}(T)$ in the direction of a field $L_{h^*}$, $h\in\g_1$. We already defined a map $\gamma$ in the projective situation :
\begin{defi}\label{gamma1}
We define $\gamma$ on $\otimes^k\gm\otimes gl(V_1,V_2)$ by 
\[\begin{array}{r}\gamma(h)(x_1\otimes\cdots\otimes x_k\otimes l)=-\sum_{i=1}^k
x_1\otimes\cdots(i)\cdots\otimes x_k\otimes (l\circ \rho_{1_*}([h,x_i]))\\
+\sum_{i=1}^k\sum_{j>i}x_1\otimes\cdots(i)\cdots\otimes\underbrace{[[h,x_i],x_j]}_{(j)}\otimes\cdots\otimes
x_k\otimes l.
\end{array}\]
for every $x_1,\cdots, x_k\in\gm$, $l\in V_1^*\otimes V_2$ and $h\in \g_1$. Then we extend it to $C^\infty(P,\otimes^k\gm\otimes gl(V_1,V_2))$ by $C^\infty(P)$-linearity.
\end{defi}
The main property of this map is the following.
\begin{prop}\label{propgamma}
For every $T\in C^\infty(P,\otimes^k\gm\otimes gl(V_1,V_2)$, one has
\[L_{h^*}\circ Q_{\omega}(T)-Q_{\omega}(T)\circ L_{h^*}=Q_\omega(L_{h^*}T+\gamma(h)T)\]
on $G_0$-equivariant functions and for every $h\in\g_1$.
\end{prop}
\begin{proof}The proof is straightforward an is similar to the corresponding one in \cite{MR1}.
\end{proof}
Moreover, we can write $\gamma$ in a very compact way. To this aim, we define a new representation of $G_0$ and $\g_0$ on $\otimes^k\gm\otimes gl(V_1,V_2)$ :
\begin{defi}
The representation $\rho_r$ of $G_0$ on $\otimes^k\gm\otimes gl(V_1,V_2)$ is defined by
\[\rho_r(a)(x_1\otimes\cdots\otimes x_k\otimes l)=Ad(a)x_1\otimes\cdots\otimes Ad(a)x_k\otimes (l\circ \rho_1(a^{-1})).\]
Note that the adjoint action of $G_0$ on $\gm$ identifies with the natural action of $CO(p,q)$.

The corresponding representation of $\g_0$ on the same space is 
\[\begin{array}{lll}\rho_{r_*}(A)(x_1\otimes\cdots\otimes x_k\otimes l)&=&\sum_{i=1}^kx_1\otimes\cdots\otimes [A,x_i]\otimes\cdots\otimes x_k\otimes l\\&&-x_1\otimes\cdots\otimes x_k\otimes l\circ \rho_{1_*}(A).\end{array}\]
for every $A\in\g_0$, $x_1,\ldots,x_k\in\gm$ and $l\in gl(V_1,V_2)$.
\end{defi}
Then we have the following immediate result
\begin{prop}\label{formrecu}
There holds
\[\gamma(h)(x_1\otimes\cdots\otimes x_k\otimes l)=\sum_{i=1}^k x_1\otimes \cdots\otimes x_{i-1}\otimes \rho_{r_*}([h,x_i])(x_{i+1}\otimes\cdots\otimes x_k\otimes l).\]
\end{prop}
Now, we are interested in the commutation relations of $\gamma$ and the representations $\rho$ and $\rho_r$. 
\begin{prop}\label{eqalg}
The following holds on $C^\infty(P,\otimes^k\gm\otimes gl(V_1,V_2))$, for $\alpha\in\{\rho,\rho_r\}$ :
\[\alpha(a)\circ \gamma(h)=\gamma(Ad(a)h)\circ \alpha(a),\]
and
\[\alpha_*(A)\circ\gamma(h)=\gamma([A,h])+\gamma(h)\circ \alpha_*(A)\]
for all $a\in G_0$, $A\in \g_0$ and all $h\in \g_1$.
\end{prop}
\begin{proof}
First notice that, since all operators under consideration are $C^\infty(P)$-linear, we only have to prove that the desired relations hold on $\otimes^k\gm\otimes gl(V_1,V_2)$. 

In order to link the relations for $\rho$ and $\rho_r$, we extend the representation $\rho_2$ of $G_0$ (see section \ref{Diffops}) to $\otimes^k\gm\otimes gl(V_1,V_2)$ in a natural way by setting
\[\rho_2(a)(x_1\otimes\cdots\otimes x_k\otimes L)=x_1\otimes\cdots\otimes x_k\otimes \rho_2(a)\circ L.\]
It is then obvious that the operators $\rho_2(a)$ and $\rho_r(b)$ commute for all $a$ and $b$ in $G_0$ and that $\rho(a)=\rho_2(a)\circ\rho_r(a).$ Therefore, we directly get
\[\alpha(a)\circ\rho_r(\exp tB)=\rho_r(a\,\exp tB\,a^{-1})\circ\alpha(a)\]
for every $a\in G_0$ and $B\in \g_0$ and $t\in\R$.
Differentiating this expression, we obtain
\[\alpha(a)\circ\rho_{r_*}(B)=\rho_{r_*}(Ad(a)B)\circ\alpha(a).\]
Now we proceed by induction on $k$. For $k=0$, the result is obvious. 
Then we set $T_1=x_1\otimes\cdots\otimes x_k\otimes l$ and $T=x_0\otimes T_1$ for $x_0,\ldots,x_k\in \gm$ and $l\in gl(V_1,V_2)$. We use Proposition \ref{formrecu} to obtain
\begin{equation}\label{eqrecu}\gamma(h)T=\rho_{r_*}([h,x_{0}])T_{1}+x_{0}\otimes\gamma(h)T_1,\end{equation}
and we have
\[\begin{array}{lll}
\alpha(a)\circ\gamma(h)T&=&\alpha(a)\circ\gamma(h)(x_0\otimes T_1)\\
&=&\alpha(a)(\rho_{r_*}([h,x_0])T_1+x_0\otimes\gamma(h)T_1)\\
&=&\rho_{r_*}(Ad(a)[h,x_0])\alpha(a)T_1+\alpha(a)(x_0\otimes\gamma(h)T_1).
\end{array}\]
The last term is equal to
\[Ad(a)x_0\otimes \alpha(a)\circ\gamma(h)T_1,\]
that is to
\[Ad(a)x_0\otimes \gamma(Ad(a)h)\circ\alpha(a)T_1\]
by induction. We then have
\[\begin{array}{lll}
\alpha(a)\circ\gamma(h)T&=&\rho_{r_*}([Ad(a)h,Ad(a)x_0])\alpha(a)T_1+Ad(a)x_0\otimes \gamma(Ad(a)h)\circ\alpha(a)T_1\\
&=&\gamma(Ad(a)h)\alpha(a)T,
\end{array}\]
by using (\ref{eqrecu}) and noticing that $\alpha(a)T=Ad(a)x_0\otimes \alpha(a)T_1$. 
\end{proof}
\begin{prop}\label{commutegamma}
There holds
\[[\gamma(h),\gamma(h')]=0\]
on $\otimes^k\gm\otimes gl(V_1,V_2)$ for all $k\in\mathbb{N}$ and $h,h'\in \g_1$.
\end{prop}
\begin{proof}
For $k=0$ or $k=1$, there is nothing to prove, since $\gamma(h)$ lowers the degree of tensors. Let us now proceed by induction.
We use the notation of Proposition \ref{eqalg}. One has then, iterating (\ref{eqrecu}), 
\begin{equation}\label{eqgamma}\gamma(h')\gamma(h)T=\gamma(h')\rho_{r_*}([h,x_0])T_1+\rho_{r_*}([h',x_0])\gamma(h)T_{1}+x_{0}\otimes\gamma(h')\gamma(h)T_{1}.\end{equation}
Using Proposition \ref{eqalg}, we obtain
\[\gamma(h')\rho_{r_*}([h,x_{0}])T_{1}=\rho_{r_*}([h,x_0])\gamma(h')T_1-\gamma([[h,x_{0}],h'])T_1.\]
It is then obvious that (\ref{eqgamma}) is symmetric in $h$ and $h'$ by induction.
\end{proof}
%%%%%
\section{Casimir-like operators}
In the papers dealing with equivariant quantization (see for instance \cite{DLO,BM,BHMP}, the existence of quantizations over vector spaces was ruled by the properties of some Casimir operators associated to the equivariance algebra. In this section we will generalize these operators to our setting. Unfortunately, we have to define them by analogy and not as true Casimir operators. Therefore, we will have to check their properties by direct computations. Hopefully, these computations are quite nice. Let us begin by what we call the flat Casimir operator. We use the basis of the algebra $\g$ defined by Proposition \ref{basis} in section \ref{algebra}.
\begin{defi}
The operator $\Cf$ is defined on $C^\infty(P,\otimes^k\gm\otimes gl(V_1,V_2))$ by 
\[\Cf=-\frac{1}{2}\rho_*(\euler)+\frac{1}{2d}\rho_*(\euler)^2+\sum_{j=1}^{\dim \h_0}\rho_*(A_j)\rho_*(A_j^*).\]
\end{defi}
The main property of this operator is the following.
\begin{prop}\label{diacasim}
The Casimir operator $\Cf$ is semi-simple. More precisely, the vector space $\otimes^k\gm\otimes gl(V_1,V_2)$ can be decomposed as an $\h_0$-representation (see section \ref{algebra}) into irreducible components (since $\h_0$ is semi-simple):
\[\otimes^k\gm\otimes gl(V_1,V_2)=\oplus_{s=1}^{n_k}I_{k,s}.\]
The restriction of $\Cf$ to $C^\infty(P,I_{k,s})$ is then a scalar multiple of the identity.
\end{prop}
\begin{proof}
The proof goes as in \cite{BM} and \cite{MR1}. Just notice that $\Cf$ is $C^\infty(P)$-linear. Thus, we only have to compute it on $I_{k,s}$ for every $s$. Then, it is easy to see that the operator $\rho_*(\euler)$ is a scalar multiple of the identity, by using the definition of $\euler$ and of $\rho$. It was proved in \cite{BM} that the last term in the expression of $\Cf$ is a scalar multiple of the Casimir operator of $\h_0$, if $\h_0$ is absolutely simple. Finally, the restriction of the Casimir operator of $\h_0$ to every irreducible representation is a scalar multiple of the identity, by Schur's Lemma. 
\end{proof}
Using the same basis, we define another operator.
\begin{defi}The operator $N^{\omega}$ is defined on $C^\infty(P,\otimes^k\gm\otimes gl(V_1,V_2))$ by
\[N^{\omega}=-2\sum_{i=1}^d\gamma(\varepsilon^i)L_{\omega^{-1}}(e_i),\]and we set
\[C^\omega=\Cf+N^\omega.\]
\end{defi}
The operator $N^{\omega}$ has an important property of invariance :
\begin{prop}\label{Goinv}
The operator $N^{\omega}$ preserves the $G_0$-equivariance of functions.
\end{prop}
\begin{proof}
The proof is exactly the same as in \cite{MR1}. This property is a consequence of the proposition \ref{eqalg} and of the fact that the invariant differentiation
preserves the $G_{0}$-equivariance. 
\end{proof}
We now have a technical lemma.
\begin{lem}
There holds on $C^\infty(P,\otimes^k\gm\otimes gl(V_1,V_2))$
\[\Cf\circ \gamma(h)-\gamma(h)\circ \Cf=2\sum_{i=1}^d\gamma(\varepsilon^i)\rho_*([h,e_i])\]
for all $h$ in $\g_1$.
\end{lem}
\begin{proof}
Using the definition of $\Cf$ and Proposition \ref{eqalg}, we directly obtain the relation
\[\begin{array}{lll}\Cf\gamma(h)&=&\gamma(h)\Cf+\gamma((-\frac{1}{2}ad(\euler)+\frac{1}{2d}ad(\euler)^2+\sum_{j=1}^{\dim \h_0}ad(A_j)ad(A_j^*))h)\\
&&+\frac{2}{2d}\gamma([\euler,h])\rho_*(\euler)+\sum_{j=1}^{\dim \h_0}(\gamma([A_j,h])\rho_*(A_j^*)+\gamma([A_j^*,h])\rho_*(A_j)).
\end{array}\]
On the other hand, using Proposition \ref{rel}, we may compute 
\[[\varepsilon^r, e_i]=\frac{1}{2d}\delta_i^r\euler-\frac{1}{2}\sum_{j=1}^{\dim h_0}(a_{ji}^rA_j^*+a_{ji}^{*r}A_j)\]
so that
\[\begin{array}{lll}2\sum_{i=1}^d\gamma(\varepsilon^i)\rho_*([h,e_i])&=&2\sum_{i,j}h_r\gamma(\varepsilon^i)\rho_*([\varepsilon^r,e_i])\\
&=&\frac{2}{2d}\gamma(h)\rho_*(\euler)-\sum_{i,j}h_r(a_{ji}^r\gamma(\varepsilon^i)\rho_*(A_j^*)+a_{ji}^{*r}\gamma(\varepsilon^i)\rho_*(A_j))\\
&=&\frac{2}{2d}\gamma(h)\rho_*(\euler)+\sum_{j=1}^{\dim h_0}(\gamma([A_j,h]\rho_*(A_j^*)+\gamma([A_j^*,h]\rho_*(A_j)), 
\end{array}\]
by using again Proposition \ref{rel}.

To sum up, we now have proved the relation
\[\begin{array}{lll}\Cf\gamma(h)&=&\gamma(h)\Cf+2\sum_{i=1}^d\gamma(\varepsilon^i)\rho_*([h,e_i])\\
&&+\gamma((-\frac{1}{2}ad(\euler)+\frac{1}{2d}ad(\euler)^2+\sum_{j=1}^{\dim \h_0}ad(A_j)ad(A_j^*))h).
\end{array}\]
We now prove that the last term vanishes by looking at the Casimir operator of the adjoint action of $\g$ :
It is given by
\[C_{ad}=\sum_{i=1}^d(ad(e_i)ad(\varepsilon^i)+ad(\varepsilon^i)ad(e_i))+\frac{1}{2d}ad(\euler)^2+\sum_{j=1}^{\dim \h_0}ad(A_j)ad(A_j^*)\]
Since $\g_1$ is an abelian subalgebra, the restriction of this operator to $\g_1$ is
\[\begin{array}{lll}C_{ad}|_{\g_1}&=&\sum_{i=1}^dad([\varepsilon^i,e_i])+\frac{1}{2d}ad(\euler)^2+\sum_{j=1}^{\dim \h_0}ad(A_j)ad(A_j^*)\\
&=&\frac{1}{2}ad(\euler)+\frac{1}{2d}ad(\euler)^2+\sum_{j=1}^{\dim \h_0}ad(A_j)ad(A_j^*).\end{array}\]
Hence, we just need to prove  
\[C_{ad}|_{\g_1}-ad(\euler)|_{\g_1}=0,\]
i.e., $C_{ad}|_{\g_1}=Id_{\g_1}$. We compute $C_{ad}$ on $\g$ in the following way :
on the one hand, since $\euler$ is in the centre of $\g_0$, we have
\[C_{ad}\euler=\sum_{i=1}^d(ad(e_i)ad(\varepsilon^i)+ad(\varepsilon^i)ad(e_i))\euler=2\sum_{i=1}^d[\varepsilon^i,e_i]=\euler,\]
by the definition of $\euler$ and Proposition \ref{basis}. On the other hand, since $\g$ is simple, we may apply Schur's Lemma. We obtain $C_{ad}=Id$ and the result follows.
\end{proof}
Now, we can come to the main result about the second Casimir operator
\begin{prop}\label{com}
There holds
\[[C^{\omega},L_{h^*}+\gamma(h)]=0\]
on $C^\infty(P,\otimes^k\gm\otimes gl(V_1,V_2))_{G_0}$, for all $h\in\g_1$.
\end{prop}
\begin{proof}
By the very definition of $C^\omega$, we have
\[\begin{array}{lll}[C^{\omega},L_{h^*}+\gamma(h)]&=&[\Cf-2\sum_{i=1}^d\gamma(\varepsilon^i)L_{\omega^{-1}}(e_i),L_{h^*}+\gamma(h)]\\
&=&[\Cf,L_{h^*}]+[\Cf,\gamma(h)]+[-2\sum_{i=1}^d\gamma(\varepsilon^i)L_{\omega^{-1}}(e_i),L_{h^*}]\\
&&+[-2\sum_{i=1}^d\gamma(\varepsilon^i)L_{\omega^{-1}}(e_i),\gamma(h)].\end{array}\]
The first and last term vanish in view of the definition of $\Cf$, $\gamma$ and of Proposition \ref{commutegamma}. We already computed the second one, and the third one is obviously equal to
\[-2\sum_{i=1}^d\gamma(\varepsilon^i)L_{[e_i,h]^*}=-2\sum_{i=1}^d\gamma(\varepsilon^i)\rho_*([h,e_i]),\]
hence the result.
\end{proof}
%%%%
\section{Construction of the quantization}\label{curvconst}
The construction is based on the eigenvalue problem for the Casimir-like operators $\Cf$ and $C^\omega$. The construction was given in \cite{MR1} in the projective case, and based on the original computations of \cite{DLO}. Actually, this construction applies to our setting. The main point is that we modified the definitions of $Q_\omega$, $\gamma$ $\Cf$ and $C^\omega$ so that Propositions \ref{propgamma}, \ref{Goinv} and \ref{com} hold true and $\Cf$ is semi-simple. We recall here the key results of the construction and we refer the reader to \cite{MR1} for the proofs. 

Recall that $\otimes^k\g_{-1}\otimes gl(V_1,V_2)$ is decomposed as a representation of $\h_0$ as the direct sum of irreducible components $I_{k,s}$. Denote by $E_{k,s}$ the space $C^\infty(P,I_{k,s})$ and by $\alpha_{k,s}$ the eigenvalue of $\Cf$ restricted to $E_{k,s}$.

As in \cite{BM,DLO}, the tree-like susbspace $ \tree(I_{k,s})$ associated to $I_{k,s}$ is defined by
\[
\tree(I_{k,s})=\bigoplus_{l\in\N}\tree^l(I_{k,s}),
\]
where $\tree^0(I_{k,s})=I_{k,s}$ and
$\tree^{l+1}(I_{k,s})=\gamma(\g_1)(\tree^l(I_{k,s}))$,
for all $l\in\N$.
The space $\tree^l(E_{k,s})$ is then defined in the same way. Since $\gamma$ is $C^\infty(P)$-linear, this space is equal to $C^{\infty}(P,\tree^l(I_{k,s}))$.

The following definition is a direct generalization of the ones of \cite{BM,
  DLO} :
\begin{defi}
An ordered pair of representations $(V_1,V_2)$
is \emph{critical} 
if there exists $k,s$ such that the eigenvalue $\alpha_{k,s}$ corresponding to an irreducible component $I_{k,s}$ of $S^k_{V_1,V_2}$
belongs to the spectrum of the restriction of $C^{\flat}$ to 
$\bigoplus_{l\geq 1}\tree^l(E_{k,s})$.
\end{defi}
We can now analyze the eigenvalue problem for
the operator $C^{\omega}$.
\begin{thm}\label{hat}
If the pair $(V_1,V_2)$ is not critical, for every $T$
in $C^{\infty}(P,I_{k,s})$, (where $I_{k,s}$ is an irreducible component of $S^k_{V_1,V_2}$) there exists a unique function $\hat{T}$ in $ C^{\infty}(P,\tree(I_{k,s}))$
 such that
\begin{equation}\label{curvP}\left\{\begin{array}{lll}
\hat{T}&=&T_k+\cdots+T_0,\quad T_k=T\\
C^{\omega}(\hat{T})&=&\alpha_{k,s}\hat{T}.
\end{array}\right.
\end{equation} 
Moreover, if $T$ is $G_0$-invariant, then $\hat{T}$ is $G_0$-invariant.
\end{thm}
%%%%
This result allows to define the main ingredient in order to define the
quantization : The ``modification map'', acting on symbols.
\begin{defi}
Suppose that the pair $(V_1,V_2)$ is not critical. Then
the map 
\[R: \oplus_{k=0}^\infty C^{\infty}(P,S^k_{V_1,V_2})\to \oplus_{k=0}^\infty C^{\infty}(P,\otimes^k\gm\otimes gl(V_1,V_2))\]
is the linear extension of the association $T\mapsto \hat{T}$.
\end{defi}
The map $M$ has the following nice property :
\begin{prop}
There holds
\begin{equation}
\label{Q}(L_{h^*}+\gamma(h))R(T)=R(L_{h^*}T),
\end{equation}
for every $h\in\g_1$ every $T\in C^{\infty}(P,S^k_{V_1,V_2})_{G_0}$ and every $k\in\N$.
\end{prop}
And finally, the main result :
\begin{thm}
If the pair $(V_1,V_2)$ is not critical, then the formula
\[Q_M: (\nabla,T)\mapsto Q_M(\nabla,T)(f)=(p^*)^{-1}[Q_{\omega}(R(p^*T))(p^*f)],\]
(where $Q_{\omega}$ is given by (\ref{Qom})) defines a natural and conformally invariant quantization.
\end{thm}
\subsection{Final remarks}
Throughout the computations, we did not use explicitly the bracket of the algebra $so(p+1,q+1)$, we only used the $|1|$-grading of this algebra and the subsequent properties. Another ingredient is the existence of a Cartan bundle associated to the $G_0$-bundle $P_0$, and of a normal Cartan connection to this bundle. Therefore, our construction can be generalized to the construction of an invariant quantization, once these data are given. 

Finally, we did not address the uniqueness problem of the quantization. But it was proved by F. Radoux that, even in the projective case, that is the most simple case, the quantization is not unique in general, due to the presence of the Weyl curvature tensor. Therefore, we conjecture that the quantization is not unique in general. It would be interesting to find a natural condition to impose to the quantization procedure in order to obtain the uniqueness that was one of the main features of the equivariant quantization problem.
%%%%
\section{Acknowledgements}
%%%%
It is a pleasure to thank C. Duval and V. Ovsienko for their fruitful comments.

F. Radoux thanks the Belgian FNRS for his Research Fellowship.


\begin{thebibliography}{10}

\bibitem{BHMP}
F.~Boniver, S.~Hansoul, P.~Mathonet, and N.~Poncin.
\newblock Equivariant symbol calculus for differential operators acting on
  forms.
\newblock {\em Lett. Math. Phys.}, 62(3):219--232, 2002.

\bibitem{BMMax}
F.~Boniver and P.~Mathonet.
\newblock Maximal subalgebras of vector fields for equivariant quantizations.
\newblock {\em J. Math. Phys.}, 42(2):582--589, 2001.

\bibitem{BM}
F.~Boniver and P.~Mathonet.
\newblock I{FFT}-equivariant quantizations.
\newblock {\em J. Geom. Phys.}, 56(4):712--730, 2006.

\bibitem{Bor}
M.~Bordemann.
\newblock Sur l'existence d'une prescription d'ordre naturelle projectivement
  invariante.
\newblock {\em Submitted for publication, math.DG/0208171}.

\bibitem{Bou1}
Sofiane Bouarroudj.
\newblock Projectively equivariant quantization map.
\newblock {\em Lett. Math. Phys.}, 51(4):265--274, 2000.

\bibitem{Bou2}
Sofiane Bouarroudj.
\newblock Formula for the projectively invariant quantization on degree three.
\newblock {\em C. R. Acad. Sci. Paris S\'er. I Math.}, 333(4):343--346, 2001.

\bibitem{Capinv}
A.~{\v{C}}ap, J.~Slov{\'a}k, and V.~Sou{\v{c}}ek.
\newblock Invariant operators on manifolds with almost {H}ermitian symmetric
  structures. {I}. {I}nvariant differentiation.
\newblock {\em Acta Math. Univ. Comenian. (N.S.)}, 66(1):33--69, 1997.

\bibitem{Capcart}
A.~{\v{C}}ap, J.~Slov{\'a}k, and V.~Sou{\v{c}}ek.
\newblock Invariant operators on manifolds with almost {H}ermitian symmetric
  structures. {II}. {N}ormal {C}artan connections.
\newblock {\em Acta Math. Univ. Comenian. (N.S.)}, 66(2):203--220, 1997.

\bibitem{CapGov}
Andreas {\v{C}}ap and A.~Rod Gover.
\newblock Tractor bundles for irreducible parabolic geometries.
\newblock In {\em Global analysis and harmonic analysis ({M}arseille-{L}uminy,
  1999)}, volume~4 of {\em S\'emin. Congr.}, pages 129--154. Soc. Math. France,
  Paris, 2000.

\bibitem{DLO}
C.~Duval, P.~Lecomte, and V.~Ovsienko.
\newblock Conformally equivariant quantization: existence and uniqueness.
\newblock {\em Ann. Inst. Fourier (Grenoble)}, 49(6):1999--2029, 1999.

\bibitem{DO1}
C.~Duval and V.~Ovsienko.
\newblock Conformally equivariant quantum {H}amiltonians.
\newblock {\em Selecta Math. (N.S.)}, 7(3):291--320, 2001.

\bibitem{DO}
C.~Duval and V.~Ovsienko.
\newblock Projectively equivariant quantization and symbol calculus:
  noncommutative hypergeometric functions.
\newblock {\em Lett. Math. Phys.}, 57(1):61--67, 2001.

\bibitem{Fox}
Daniel J.~F. Fox.
\newblock Projectively invariant star products.
\newblock {\em IMRP Int. Math. Res. Pap.}, (9):461--510, 2005.

\bibitem{Hansoul}
Sarah Hansoul.
\newblock Projectively equivariant quantization for differential operators
  acting on forms.
\newblock {\em Lett. Math. Phys.}, 70(2):141--153, 2004.

\bibitem{sarah}
Sarah Hansoul.
\newblock Existence of natural and projectively equivariant quantizations.
\newblock {\em Adv. Math.}, 214(2):832--864, 2007.

\bibitem{Kobabook}
Shoshichi Kobayashi.
\newblock {\em Transformation groups in differential geometry}.
\newblock Springer-Verlag, New York, 1972.
\newblock Ergebnisse der Mathematik und ihrer Grenzgebiete, Band 70.

\bibitem{Koba1}
Shoshichi Kobayashi and Tadashi Nagano.
\newblock On filtered {L}ie algebras and geometric structures. {I}.
\newblock {\em J. Math. Mech.}, 13:875--907, 1964.

\bibitem{LO}
P.~B.~A. Lecomte and V.~Yu. Ovsienko.
\newblock Projectively equivariant symbol calculus.
\newblock {\em Lett. Math. Phys.}, 49(3):173--196, 1999.

\bibitem{Leconj}
Pierre B.~A. Lecomte.
\newblock Towards projectively equivariant quantization.
\newblock {\em Progr. Theoret. Phys. Suppl.}, (144):125--132, 2001.
\newblock Noncommutative geometry and string theory (Yokohama, 2001).

\bibitem{Loubon}
S.~E. Loubon~Djounga.
\newblock Conformally invariant quantization at order three.
\newblock {\em Lett. Math. Phys.}, 64(3):203--212, 2003.

\bibitem{MR1}
P.~Mathonet and F.~Radoux.
\newblock Cartan connections and natural and projectively equivariant
  quantizations.
\newblock {\em J. Lond. Math. Soc. (2)}, 76(1):87--104, 2007.

\bibitem{MR2}
P.~Mathonet and F.~Radoux.
\newblock On natural and conformally equivariant quantizations.
\newblock 2008.

\bibitem{MR}
Pierre Mathonet and Fabian Radoux.
\newblock Natural and projectively equivariant quantiations by means of cartan
  connections.
\newblock {\em Lett. Math. Phys.}, 72(3):183--196, 2005.

\end{thebibliography}
\end{document}